\documentclass[13pt, reqno]{amsart}

\usepackage{graphicx}

\usepackage{mathrsfs}
\usepackage{amssymb}
\usepackage{amsfonts, mathtools}
\usepackage{amssymb,amsmath,amsthm,tikz-cd, tikz}

\usepackage{hyperref}
\usepackage{palatino}
\usepackage{wrapfig}
\usepackage{subcaption}

\usetikzlibrary{patterns, matrix,shapes,arrows,calc,topaths,intersections,positioning,decorations.pathreplacing,decorations.pathmorphing,fit}

\newif\ifpaper
\papertrue

\usepackage{color}

\def\name{Peng Zhou}

\newcommand{\congto}{\xto{\sim}}
\newcommand*{\hrlen}{5}
\newcommand*{\hramp}{3}

\tikzset{
asdstyle/.style={blue,thick},
righthairs/.style={postaction={decorate,draw,decoration={border,amplitude=\hramp,segment length=\hrlen,angle=-90,pre=moveto,pre length=\hrlen/2}}},
lefthairs/.style={postaction={decorate,draw,decoration={border,amplitude=\hramp,segment length=\hrlen,angle=90,pre=moveto,pre length=\hrlen/2}}},
righthairsnogap/.style={postaction={decorate,draw,decoration={border,amplitude=\hramp,segment length=\hrlen,angle=-90}}},
lefthairsnogap/.style={postaction={decorate,draw,decoration={border,amplitude=\hramp,segment length=\hrlen,angle=90}}},
graphstyle/.style={thick},
arrowstyle/.style={thick,decorate,decoration={snake,amplitude=1.7,segment length=10pt,post length=.5mm,pre length=0}},
genmapstyle/.style={thick,-stealth'},
arrhdstyle/.style={thick},
exceptarcstyle/.style={red, ultra thick},
dualquiverstyle/.style={thick,->},
patstyle/.style={pattern color = gray, pattern = north east lines, opacity=0.3}
}

\title[Twisted Polytopes and CCC] {Twisted Polytope Sheaves and Coherent-Constructible Correspondence for Toric Varieties}
\author{\name}
\date{\today}

\DeclareMathOperator{\Hom}{Hom}

\renewcommand{\hom}{\Hom}

\DeclareMathOperator{\cone}{cone}

\DeclareMathOperator{\Spec}{Spec}

\DeclareMathOperator{\supp}{Supp}

\DeclareMathOperator{\colim}{colim}

\DeclareMathOperator{\Lag}{Lag}
\DeclareMathOperator{\Leg}{Leg}
\DeclareMathOperator{\Lie}{Lie}
\DeclareMathOperator{\End}{End}
\DeclareMathOperator{\Vect}{Vect}

\newcommand{\tforall}{\quad \text{ for all } \;}
\newcommand{\orn}{\mathfrak{or}}

\newcommand{\pa}{\partial}

\newcommand{\la}{\langle}
\newcommand{\ra}{\rangle}
\newcommand{\bu}{\bullet}
\newcommand{\ot}{\otimes}

\newcommand{\uhom}{{\mathscr{H}om}}

\renewcommand{\ss}{\subsection}
\renewcommand{\cong}{\simeq}

\newcommand{\ccal}{\mathcal{C}}
\newcommand{\dcal}{\mathcal{D}}

\newcommand{\lcal}{\mathcal{L}}

\newcommand{\ocal}{\mathcal{O}}
\newcommand{\pcal}{\mathcal{P}}

\newcommand{\scal}{\mathcal{S}}
\newcommand{\tcal}{\mathcal{T}}

\newcommand{\wb}{\overline}

\newcommand{\C}{\mathbb{C}}
\newcommand{\D}{\mathbb{D}}

\newcommand{\R}{\mathbb{R}}

\newcommand{\Z}{\mathbb{Z}}

\newcommand{\F}{\mathbb{F}}

\renewcommand{\P}{\mathbb{P}}

\newcommand{\xto}{\xrightarrow}

\newcommand{\lra}{\leftrightarrow}
\newcommand{\LRA}{\Leftrightarrow}

\newcommand{\RM}{\backslash}
\newcommand{\into}{\hookrightarrow}

\newcommand{\bea}{\begin{eqnarray*} }
\newcommand{\eea}{\end{eqnarray*} }
\newcommand{\be}{\begin{equation} }
\newcommand{\ee}{\end{equation} }
\newcommand{\bp}{\begin{proposition}}
\newcommand{\ep}{\end{proposition}}
\newcommand{\bt}{\begin{maintheo}}
\newcommand{\et}{\end{maintheo}}
\newcommand{\bpf}{\begin{proof}}
\newcommand{\epf}{\end{proof}}
\newcommand{\bl}{\begin{lemma}}
\newcommand{\el}{\end{lemma}}
\newcommand{\bc}{\begin{corollary}}
\newcommand{\ec}{\end{corollary}}
\newcommand{\bd}{\begin{definition}}
\newcommand{\ed}{\end{definition}}
\newcommand{\bee}{\begin{eqnarray} }
\newcommand{\eee}{\end{eqnarray} }
\newcommand{\brem}{\begin{remark}}
\newcommand{\erem}{\end{remark}}
\newcommand{\bex}{\begin{example}}
\newcommand{\eex}{\end{example}}

\newcommand{\bma}{\begin{bmatrix}}
\newcommand{\ema}{\end{bmatrix}}
\newcommand{\bcs}{\begin{cases}}
\newcommand{\ecs}{\end{cases}}
\newcommand{\bcd}{\begin{tikzcd}}
\newcommand{\ecd}{\end{tikzcd}}

\newtheorem{maintheo}{Theorem}
\newtheorem{theo}{Theorem}[section]
\newtheorem{lemma}[theo]{Lemma}
\newtheorem{corollary}[theo]{Corollary}
\newtheorem{proposition}[theo]{Proposition}

\theoremstyle{definition}
\newtheorem{definition}[theo]{Definition}
\newtheorem{example}[theo]{Example}
\newtheorem{remark}[theo]{Remark}

\theoremstyle{plain}

\begin{document}

\begin{abstract}
Given a smooth projective toric variety $X_\Sigma$ of complex dimension $n$,  Fang-Liu-Treumann-Zaslow \cite{FLTZ} showed that  there is a quasi-embedding of the differential graded (dg) derived category of coherent sheaves $Coh(X_\Sigma)$ into the dg derived category of constructible sheaves on a torus $Sh(T^n, \Lambda_\Sigma)$. Recently, Kuwagaki \cite{Ku2} proved that the quasi-embedding is a quasi-equivalence, and generalized the result to toric stacks. Here we give a different proof in the smooth projective case,  using non-characteristic deformation of sheaves to find twisted polytope sheaves that co-represent the stalk functors. 
\end{abstract}

\maketitle

\section{Introduction}
Toric varieties are certain compactifications of the complex torus $(\C^*)^n$. They provide many interesting examples, and can be studied in various ways, using algebraic geometry, symplectic geometry or combinatorics. 

\begin{figure}[h]
\centering
%
%
%
%
%
%
%
%
%
%
%
%

\includegraphics[width=\textwidth]{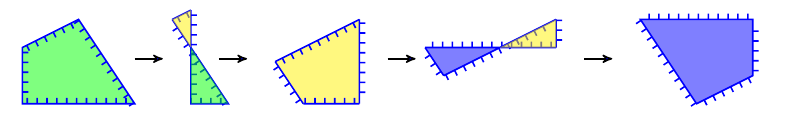}

\caption{\label{fig:twist} Twistings of a convex polytope. As one pushes the edges of the polytopes, a certain edge will shrink to zero-length then reappear on the other side. Note the co-directions of the edges, indicated by hairs, remain fixed. The corresponding twisted polytope sheaves have stalks on green, yellow and blue regions as $\C, \C[1], \C[2]$, respectively.}
\end{figure}

For example, let $X_\Sigma$ be a smooth projective toric vareity corresponding to a fan $\Sigma$, and $L$ an ample line bundle
with a lifting of the $(\C^*)^n$-action.  Then there is a convex polytope $\Delta_L$ in $\R^n$, where $\R^n$ is identified with the dual Lie algebra $\Lie(T^n)^\vee$ and $T^n = (U(1))^n$ is the maximal compact real subgroup of $(\C^*)^n$. The convex polytope $\Delta_L$ can be understood in the following ways, 
\begin{enumerate}
\item Algebraically, $\Delta_L$ is the convex hull of the characters appearing in the weight decomposition of $H^0(X, L)$ under the $(\C^*)^n$-action. 
\item Symplectically, $\Delta_L$ is the moment polytope of the Hamiltonian action $T^n$ on $(X, \omega)$, where $\omega$ is a symplectic 2-form with $[\omega] = c_1(L)$. 
\item Combinatorially, $\Delta_L$ is the intersection of half-spaces $Q_\rho=\{x \in \R^n \mid \la x, v_\rho \ra \leq a_\rho\}$, one for each compactifying divisor $D_\rho$ of $X$ given by a vector $v_\rho \in \Z^n$, and $a_\rho$ is   the vanishing order of the invariant (meromorphic) section along the divisor $D_\rho$.
\end{enumerate}
In the case where $L$ is not ample, the corresponding polytope becomes a `twisted polytope', as explained in Figure \ref{fig:twist}. The name originates from the paper of Karshon and Tolman \cite{KT}, where they generalized the moment map to the case where $\omega$ is degenerate. 

The above correspondence between equivariant line bundles and twisted polytopes enjoys a categorification under the name of (equivariant) {\em Coherent-Constructible-Correspondence (CCC)}. 
\bt [\cite{FLTZ}] \label{t:FLTZ}
If $X$ is a proper toric variety, there is a corresponding conical Lagrangian $\Lambda_\Sigma \subset T^* \R^n$ and an equivalence of derived (or rather, triangulated dg categories)
\[ \kappa: Perf_T(X_\Sigma) \congto Sh_{cc}(\R^n, \Lambda_\Sigma) \]
where 
\begin{itemize}
\item $Perf_T(X_\Sigma)$ is the triangulated dg category of perfect complexes of torus-equivariant coherent sheaves on $X_\Sigma$. 
\item $Sh_{cc}(\R^n, \Lambda_\Sigma)$ is the triangulated dg category of constructible sheaves on $\R^n$ which are compactly supported, whose  singular supports lie in $\Lambda_\Sigma$. 
\end{itemize}
\et

The equivariant CCC implies that there is a quasi-embedding for the non-equivariant case: 
\bp[\cite{Tr},Proposition 2.4, 2.7]
Let $\pi: \R^n \to T^n \cong \R^n / \Z^n$ be the projection. Then there exists a functor $\wb \kappa$ and commutative diagrams
\[ 
\begin{tikzcd}
Perf_T(X_\Sigma) \arrow{r}{\sim\, \kappa} \arrow{d}{\text{forget}} & Sh_{cc}(\R^n, \Lambda_\Sigma) \arrow{d}{\pi_!} \\
Perf(X_\Sigma)  \arrow[hookrightarrow]{r}{\wb \kappa} & Sh(T^n, \wb \Lambda_\Sigma), 
\end{tikzcd}
\]
\ep
\brem
When $X_\Sigma$ is smooth, the homotopy category of $Perf_T(X_\Sigma)$ (resp. $Perf(X_\Sigma)$) coincide with the usual $D^bCoh_T(X_\Sigma)$ (resp. $D^b Coh(X_\Sigma)$). 
\erem

\brem
Under the quotient map $\pi: \R^n \to T^n$, all the upstairs objects in $\R^n$ are unadorned, and downstairs objects in $T^n$ have overlines. 
\erem
And it is conjectured that this quasi-embedding is a quasi-equivalence. The conjecture has been verified in certain cases by Treumann \cite{Tr}, Scherotzke-Sibilia \cite{SS} and Kuwagaki \cite{Ku1}. Recently, it has been fully proven by Kuwagaki \cite{Ku2} in the generality of toric stacks, using gluing descriptions of $\infty$-categories on both sides. 

In this paper, we prove the non-equivariant CCC for smooth projective toric varieties, by showing the $\wb \kappa$-images of line bundles generate the constructible sheaf category. 
\begin{maintheo}
Let $X_\Sigma$ be a smooth projective toric variety of complex dimension $n$, then there is an quasi-equivalence of category
\[\wb \kappa: Coh(X_\Sigma) \congto Sh(T^n, \wb \Lambda_\Sigma) \]
where $\wb \Lambda_\Sigma$ is a conical Lagrangian in $T^*T^n$. 
\end{maintheo}

The key part of the proof is as following. For any point $\theta \in T^n$, there is a constructible sheaf  $\wb P_{[\theta]}$ on $T^n$ as the $\wb \kappa$-image of a certain line bundle (c.f. Definition \ref{d:tw-div}), such that for any sheaf $ F \in Sh(T^n, \wb \Lambda_\Sigma)$, its stalk at the point $\theta$ can be computed by
\be \label{eq:key} F_\theta \cong \hom(\wb P_{[\theta]}[-n], F). \ee
This immediately implies that if  $F$ satisfies $\hom(\wb \kappa(L), F) = 0$ for all the line bundles $L$ on $X_\Sigma$, then $F=0$. In other words, the stalk functors in $Sh(T^n, \wb \Lambda_\Sigma)$ are co-represented by $\wb \kappa$-images of line bundles on $X_\Sigma$. We thank David Treumann for the suggestion of  co-representing the stalk functors. 

The quasi-isomorphism \eqref{eq:key} is due to a non-characteristic deformation argument for constructible sheaf. We define a 1-parameter family of sheaves $\{P_t\}_{t \in [0,1]}$, such that 
\begin{enumerate}
\item $P_0 = j_{B!} \C_B$, where $B$ is a small enough convex open set around $\theta$, such that $F_\theta \cong \Gamma(B, F) \cong \hom(P_0, F)$.
\item $P_1 = \wb P_{[\theta]}[-n]$. 
\item For $t \in (0,1)$, take the linear interpolation between $P_0$ and $P_1$, and show that $SS^\infty(P_t) \cap \wb \Lambda_\Sigma^\infty = \emptyset$. 
\end{enumerate}
By the non-characterstic deformation lemma\footnote{One needs to be careful about the endpoint $t=1$, since $SS^\infty(P_1) \cap \wb \Lambda_\Sigma^\infty \neq \emptyset$. The non-characteristic deformation lemma for sections over open sets, Proposition \ref{p:ks-nc}, avoids this problem. },   $\hom(P_t, F)$ is invariant during the deformation, hence we get \eqref{eq:key}.

\begin{figure}
\centering
%
\includegraphics[width=\textwidth]{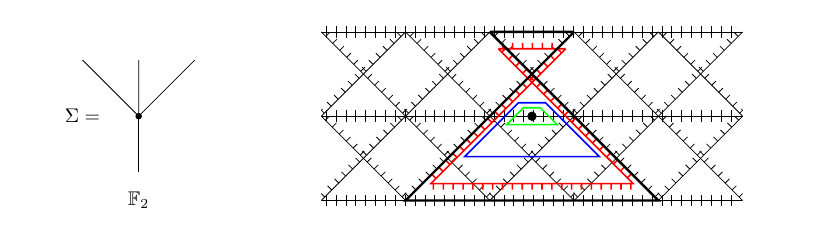}
\caption{\label{fig:expanding} Expanding family of twisted polytope sheaves on $\R^2$.}

\end{figure}

\bex [Expanding family of twisted polytope sheaves]
The example of Hirzebruch surface $\F_2$, see Figure \ref{fig:expanding}. Here we describe the sheaf $P_{[x]}$ upstairs in $\R^2$, where $\wb P_{[\pi( x)]}:= \pi_* P_{[x]}$ and $\pi:\R^2 \to T^2$ is the quotient map. The point we want to probe is at $x=(-0.5, 0)$, marked in black. The green, blue, red and black curves are the boundaries of the twisted polytopes in the interpolating family $P_t$. The green and blue ones are still open convex polytope, the red and black ones are twisted. We marked the direction of the singular support for the sheaf $P_{\text{red}}$ corresponding red curve, and note that  $SS^\infty(P_{\text{red}}) \cap \Lambda^\infty = \emptyset$. 
\eex

\brem
The collection of line bundles as the $\wb \kappa$-preimages of $\{ \wb P_{[\theta]} \}$ is a finite collection, since sheaves in $Sh(T^n, \wb \Lambda_\Sigma)$ admits a finite stratification depending only on $\wb \Lambda_\Sigma$. In the case of $\P^n$, they turn out to be $\ocal(1), \cdots, \ocal(n+1)$, and form an exceptional collection. However, for general toric variety, even smooth Fano ones, the collection of line bundles cannot always be an exceptional collection \cite{E, HP}. 
 
However, for any collection of line bundles $L_1, \cdots, L_N$, Craw-Smith \cite{CS} considered the endomorphism algebra $A = \End(\oplus_{i=1}^N L_i)$ and defined a bound quiver of sections $(Q,R)$ associated with $A$. It would be interesting to study our collection of line bundles using this quiver approach. 
\erem

The outline of the paper is the following. In section 2 we review the background on constructible sheaves, in particular state the non-characteristic deformation lemma. In section 3, we review toric geometry and state the relevant results from \cite{FLTZ, Tr} on the equivariant and non-equivariant CCC functors.  In section 4, we prove our main theorem. 

\ss{Acknowledgements}
It is a pleasure to thank Xin Jin, Linhui Shen, Lei Wu, Elden Elmanto and Dima Tamarkin for many helpful discussions, and David Treumann and David Nadler for their interests in this work. I am grateful for Elden for carefully reading the draft and giving many useful comments. I am grately indebted to my advisor Eric Zaslow for inspirations and encouragements (and patience!). The discussion with David Treumann at IAS inspired the current approach, which uses the twisted polytope sheaves to corepresent the stalk functor. 


\section{Constructible Sheaves}
Here we give a quick review on constructible sheaves, following the introduction of \cite{N} and the Appendix of \cite{STW} very closely.  For a thorough account on constructible sheaf theory and its relation with Fukaya category, see \cite{KS, S} and \cite{NZ}. 

We first give the categorical background, especially the definition for dg enhancement of the triangulated derived categories. Then we give some useful formulae for practical computations. Then we take a brief detour in symplectic geometry to define conical Lagrangian in cotangent bundles, so that we can define the singular support $SS(F)$ of a constructible sheaf $F$. Finally, we discuss various {\em non-characteristic}  deformations results, where the sections  $\Gamma(U_t, F)$, or more generally the hom-complexes $\hom(P_t, F)$, is invariant upto quasi-isomorphism as $U_t$ or $P_t$ vary. 

\ss{Classical and differential graded derived categories of sheaves}
Let $X$ be a topological space. The poset (viewed as a category) $Open(X)$ has objects of open subsets, and partial orderings (morphisms) are given by inclusions of open subsets. Let $\Vect$ be the abelian category of  $\C$-vector spaces.
\begin{itemize}
\item  A presheaf $F$ valued in $\Vect$ is a functor 
\[ F: Open(X)^{o} \to \Vect. \]
A presheaf $F$ is a sheaf, if for any collection of open subsets $\{U_i\}_{i \in I}$, we have an exact sequence
\[ 0 \to F( \bigcup_{i  } U_i) \to \prod_{i  } F(U_i) \to \prod_{i,j} F(U_i \cap U_j). \]
\item Let $C(X)$ be the abelian category of complexes of sheaves on $X$, with morphisms being degree-zero chain maps. 
\item Let $K(X)$ be the homotopy category of $C(X)$, where objects are the same as $C(X)$, ands morphisms are chain maps upto homotopy equivalences
\[ \Hom_{K(X)}(F^\bu, G^\bu) := \Hom_{C(X)}(F^\bu, G^\bu) / \sim \]
where $\varphi_1 \sim \varphi_2$ if  $\varphi_1 - \varphi_2 = d \circ h - h \circ d$ for some degree $-1$ map $h: F^\bu \to G^{\bu-1}$. 
\item  The derived category $D(X)$ is obtained from $K(X)$ by inverting quasi-isomorphisms. The bounded derived category $D^b(X)$ is defined to be the full subcategory of complexes with bounded cohomologies. 
\end{itemize}

To define constructibility, let $X$ be a real analytic manifold. We fix an algebro-geometric category $\ccal$, e.g., the category of subanalytic sets. 
A {\em Whitney stratification} $\scal = \{\scal_\alpha\}_{\alpha \in A}$ of $X$ by $\ccal$-submanifolds, is a decomposition $X = \sqcup_{\alpha \in A} \scal_\alpha$ into disjoint strata $\{\scal_\alpha\}$ indexed by $A$, where $\scal_\alpha$ are locally closed $\ccal$-submanifolds and $\scal_\alpha \cap \wb{\scal_\beta} \neq \emptyset$ if and only if $\scal_\alpha \subset \wb{\scal_\beta}$,  and any pair of distinct strata $(\scal_\alpha, \scal_\beta)$ satisfies the {\em Whitney condition}, that is, if a sequence $\{x_n \in \scal_\alpha\}$ and a sequence $\{y_n \in \scal_\beta\}$ converges to a point $y \in \scal_\beta$, such that in some local coordinate patch, the secant lines $\wb{x_i y_i}$ converges to a line $l$ and the tangent planes $T_{x_i} \scal_\alpha$ converges to a plane $\tau$, then $l \subset  \tau$. 

Let $\scal=\{\scal_\alpha\}$ be a Whitney stratification of $X$. An object $F^\bu \in D(X)$ is said to be $\scal$-constructible, if the restrictions $H^i(F^\bu)|_{\scal_\alpha}$ of its cohomology sheaves to the strata of $\scal$ are finite-rank and locally constant. We denote by $D_\scal(X)$ the full subcategory of $D(X)$ spanned by $\scal$-constructible objects, and denote by $D_c(X)$ the full subcategory of $D(X)$ spanned by  constructible sheaves for some Whitney stratification. If the stratification $\scal$ is finite, then by the finite rank condition on cohomology sheaves, $D_c(X) \subset D^b(X)$. 

Next, we define the differential graded(dg) derived category. For a review on dg category and dg quotient construction, see \cite{Ke} and \cite{Dr}. 
\begin{itemize}
\item The (naive) dg category $Sh_{naive}(X)$ has objects as chain complexes of sheaves $F^\bu$,  same as $C(X)$, and morphisms are chain complexes, with the degree $n$ element as
\[ \hom^n_{Sh_{naive}(X)}(F^\bu, G^\bu):= \prod_{i \in \Z} \hom_{X}(F^i, G^{n+i}) \]
and differentials are given by
\[ d^n:  \hom^n_{Sh_{naive}(X)}(F^\bu, G^\bu) \to  \hom^{n+1}_{Sh_{naive}(X)}(F^\bu, G^\bu), \quad \varphi \mapsto d_F \circ \varphi - (-1)^n  \varphi \circ d_G. \]
\item The dg derived category $Sh(X)$ is the dg quotient of  $Sh_{naive}(X)$ by the full subcategory spanned by acyclic objects \cite{Dr}. This is a triangulated dg category whose cohomology category $H^0(Sh(X))$ is canonically equivalent to the derived category $D(X)$ as a triangulated category. 
\item The dg derived category of bounded constructible categories $Sh_c^b(X)$ is the full subcategory of $Sh(X)$, whose objects $F^\bu$ have bounded constructible cohomology sheaves. For a fixed Whitney stratification $\scal$ of $X$, the $\scal$-constructible dg derived category $Sh_\scal(X)$ is the full dg subcategory of $Sh(X)$ spanned by objects projecting to $D_\scal(X)$. 
\item A dg functor $F: \ccal \to \dcal$ of dg categories is a quasi-embedding (resp. quasi-equivalence) if and only if the induced cohomological functor $H(F): H(\ccal) \to H(\dcal)$ is an embedding (resp. equivalence). 
\end{itemize}

In this paper, we will only work with constructible sheaves, and will omit the subscript $c$ for constructibility. To simplify notation, we use sheaf $F$ to mean complex of sheaves $F^\bu$,  $\hom_X$ to mean hom-complex $\hom^\bu_{Sh(X)}$. 

\ss{Useful Formulae for Computations}
Inspite of the abstract categorical definitions, constructible sheaves enjoy many functorial properties which faciliates actual computations. Here we give some useful formulae and examples. 

We use $f^*, f_*, f^!, f_!, \uhom, \ot$ to mean the corresponding dg derived functors: 
\bea
  - \ot F: & Sh(X) \lra Sh(X) & : \uhom(F, -)\\
  f^*:& Sh(X) \lra Sh(Y) &: f_* \\
  f_!: &Sh(Y) \lra Sh(X)&: f^! 
\eea
where $f: Y \to X$ is a map of real analytic manifolds. 

The Verdier duality $\D: Sh(X)^o \to Sh(X)$ is an anti-involution. It interchanges shriek with star
\[ \D\D = id, \quad f_! = \D f_* \D, \quad f^! = \D f^* \D. \]
The shrieks and stars are directly related in two cases: when $f$ is proper $f_! = f_*$; when $f$ is a smooth morphism of relative dimension $d_f$, $f^!(-) \cong f^*(-)\otimes  \omega_{Y/X} \cong f^*(-)\otimes \orn_{Y/X} [d_f]$, where $\orn_{Y/X}$ is the orientation sheaf of the fiber.

Given an open subset $U$ of $X$ and its closed complement $Z$,  
\[  \text{open inclusion:} \quad U \xhookrightarrow{j} X \xhookleftarrow{i} Z, \quad \text{closed inclusion},  \]
we have $j^* = j^!$ and $i_*=i_!$. Furthermore, there are exact triangles
\[ i_! i^! \to id \to j_* j^* \xto{[1]}, \quad  j_! j^! \to id \to i_* i^* \xto{[1]}.  \]
These are sheaf-theoretic incarnations of excisions: applied to the constant sheaf on $X$ and taking global sections, we get
\[ H^*(Z, i^! \C) \to H^*(X, \C) \to H^*(U, \C) \xto{[1]}, \quad H^*_c(U, \C) \to H^*_c(X, \C) \to H^*_c(Z, \C) \xto{[1]}. \]

If $Y$ is a locally closed $\ccal$-submanifold of $X$, we use $j_Y: Y \into X$ to denote the inclusion. Let $\C_Y \in Sh(Y)$ denote the constant sheaf on $Y$, and $\omega_Y = \D \C_Y$ be the Verdier dualizing complex of $Y$, then $\omega_Y$ is the canonically isomorphic to the shifted orientation sheaf $ \orn_Y[\dim Y]$  on $Y$. The standard sheaf on $Y$ is $j_{Y*} \C_Y$, and the costandard sheaf on $Y$ is $j_{Y!} \omega_Y$.

The constructible sheaves can be `constructed' by taking shifts and mapping cones of certain finite collection of sheaves. 
Let $\tcal = \{ \tau_\alpha\}$ be a triangulation of $X$ by simplices $j_\alpha: \tau_\alpha \into X$, where each $\tau_\alpha$ is the embedding image of some open simplex. 
We denote by  $\ccal_*(\tcal)$ the the full dg subcategory of $Sh_\tcal(X)$ spanned by {\em standard objects} ${j_\alpha}_* \C_{\tau_\alpha}$. The morphisms between standard objects are quasi-isomorphic to complexes concentrated at degree zero. 
\[ \hom_{Sh_\tcal(X)}({j_\beta}_* \C_{\tau_\beta},{j_\alpha}_* \C_{\tau_\alpha}) \cong \begin{cases}
 \C & \text{ if $\beta \geq \alpha$ } \\
 0  & \text{ else} 
 \end{cases}
\]
where $\alpha \leq \beta$ if $\tau_\alpha \subset \wb{\tau_\beta}$.

Applying Verdier duality to the standard sheaves, we get the costandard sheaves $\D ({j_\alpha}_* \C_{\tau_\alpha}) = {j_\alpha}_! \omega_{\tau_\alpha}$. We have
\[ \hom_{Sh_\tcal(X)}({j_\alpha}_! \omega_{\tau_\alpha}, {j_\beta}_1 \omega_{\tau_\beta}) \cong \begin{cases}
 \C & \text{ if $\alpha \leq \beta$ } \\
 0  & \text{ else.} 
 \end{cases}
\]
We denote by  $\ccal_!(\tcal)$ the the full dg subcategory of $Sh_\tcal(X)$ spanned by {\em costandard objects} ${j_\alpha}_! \omega_{\tau_\alpha}$

\bl[\cite{N}, Lemma 2.3.1]
$Sh_\tcal(X)$ is the triangulated envelope of $\ccal_*(\tcal)$ (resp. $\ccal_!(\tcal)$). 
\el

\bex
Let $Y = (0,1) \subset \R$, then $j_{Y*} \C_Y$ is the constant sheaf with stalk $\C$ supported on the closed interval $[0,1]$, and $j_{Y!} \omega_Y$ is the costandard sheaf with stalk isomorphic to $\C[1]$ supported on the open interval $(0,1)$.
\eex
\bex  Let $f: \{0\} \into \R^n$, then 
\[ f^!(\C_{\R^n}) = \D f^* \D \C_{\R^n} =  \D f^*  (\C_{\R^n}[n]) = \D (\C_{\{0\}}[n]) = \C_{\{0\}}[-n] \] 
where we have identified the orientation sheaf on $\R^n$ with the constant sheaf by chosing an orientation on $\R^n$. 
\eex


\ss{Conical Lagrangian and Singular Support}
In this subsection, we define singular supports of constructible sheaves. Roughly speaking,  singular supports encode the `positions and directions' where  sections `fail to propagate'. We first need to introduce notations from symplectic and contact geometry. 

Let $X$ be a smooth manifold, $T^*X$ its cotangent bundle with the canonical one-form $\lambda = p d q$ and the canonical sympletic two form $\omega = d \lambda = dp \wedge dq$. Let $\dot{T}^*X = T^*X \RM X$, where $X$ is identified with the zero section in $T^*X$.
Let $T^\infty X = \dot T^*X / \R_{>0}$, where $\R_{>0}$ acts by fiberwise dilation. There is a natural fiberwise compactifiation of $T^*X$ to $\wb T^* X$, where $T^\infty X$ corresponds to the divisor at infinity $\wb T^*X \RM T^*X$ (\cite{NZ}, \S 5.1.1). 

A contact manifold $(M,\xi)$ is a smooth manifold of odd dimension $2m+1$, with a smooth rank $2m$ subbundle $\xi \subset TM$, called a {\em hyperplane distribution}, such that locally $\xi = \ker(\alpha)$ for some one-form $\alpha$ and $\alpha \wedge (d\alpha)^m \neq 0$. Such a one-form $\alpha$ is called a {\em contact form}. The {\em Reeb vector field} with respect to a contact form $\alpha$ is the unique vector field $R$ such that $\iota_R \alpha =1$ and $\iota_R d\alpha = 0$. A contactomorphism between contact manifolds is a diffeomorphism that preserves the hyperplane distributions. A Legendrian submanifold $\lcal$ of $M$ is an $m$-dimensional submanifold such that $T \lcal \subset \ker(\alpha) \cap \ker(d \alpha)$. 

The divisor  $T^\infty X$ at infinity of the compactification  $\wb T^*X$ has a natural contact structure defined in the following way:  Fix any smooth section $H$ of the $\R_{>0}$-bundle $\dot T^*X \to T^\infty X$, then $T^\infty X$ is diffeomorphic to $H$ by the section map;  the canonical one-form $\lambda$ of $T^*X$ restricts to a contact form $\alpha$ on $H$, hence induces a contact structure $\xi$ on $T^\infty X$.  If we fix a Riemmanian metric on $X$, then the section $H$ can be taken as the {\em unit cosphere bundle}
\[ S^*X = \{ (x, \eta) \in T^*X \mid \|\eta\|=1\}. \]
The Reeb flow on $S^*X$ is the unit geodesic flow. We will identify $S^* X$ and $T^\infty X$. 

\bex
The simplest example contact manifold is the 1-jet bundle on $\R^n$: $J^1 \R^n := T^*_{(x,y)} \R^n \times \R_z$, and one choice of the contact form can be taken as $\alpha = z - \sum_{i=1}^n y_i d x_i$ and the corresponding Reeb flow is $\pa_z$. 
\eex

A conical Lagrangian $\Lambda \subset T^* X$ is a Lagrangian (possibly singular) invariant under the $\R_{>0}$-action. A homogenous conical Lagrangian is a one contained in $\dot T^* X$. Given a conical Lagrangian $\Lambda$, we define the associated Legendrians as 
\[ \Leg(\Lambda) = \Lambda^\infty = (\Lambda \RM X) / \R_{>0} \subset T^\infty X. \] 
Conversely, given a Legendrian $\lcal \subset T^\infty X$, we use $\Lag(\lcal)$ to denote the homogeneous conical Lagrangian in $\dot T^*X$ as the preimage of the quotient $\dot T^* X \to T^\infty X$.

Let $\scal = \{ \scal_\alpha \}_{\alpha \in A}$ be a Whitney stratification of $X$, then there is a canonical conical Lagrangian associated to $\scal$, 
\[ \Lambda_\scal : = \bigcup_{\alpha \in A} T_{\scal_\alpha}^* X\]
where $T_N^*X= \{ (x, \eta) \in T^*X \mid x \in N, \eta|_{TN} = 0 \}$ denotes the conormal bundle of a submanifold $N \subset X$. 

Let $F$ be a $\scal$-constructible sheaf in $Sh_\scal(X)$ for a Whitney stratification $\scal$. The singular support $SS(F)$ is a (singular) conical Lagrangian contained in $\Lambda_\scal$ defined in the following way: a point $(x,\eta) \in T^*X$ is {\bf not} in the singular support $SS(F)$, if there is a small open ball $B(x,\epsilon)$ around $x$, and a Morse function $f: B(x,\epsilon) \to \R$, with $f(x)=0$ and $df(x) = \eta$, such that for any $0 < \delta \ll 1$, the canonical restriction morphism
\[ \Gamma(f^{-1}(-\infty, \delta), F) \to \Gamma(f^{-1}(-\infty, -\delta), F)  \]
is a quasi-isomorphism. We use 
\[ SS^\infty(F) := (SS(F))^\infty = (SS(F) \RM X) / \R_{>0} \]
 to denote the Legendrian in $T^\infty X$ associated to the conical Lagrangian $SS(F)$ in $T^* X$. 
\bex
Let $j: U=B(0,1) \into \R^2$ be the inclusion of an open unit ball in $\R^2$. Then $j_* \C_U$ is supported on the closed set $\wb U$, with singular support at infinity as 
\[ SS^\infty(j_* \C_U) = \{ (x, \eta) \in S^*\R^2 \mid x \in \pa U, \eta = -d |x| \} = \quad \tikz[baseline=-3pt] {\draw [lefthairs] (0,0) circle (0.5);}. \]
And $j_! \C_U$ is supported on the open set $U$, with singular support at infinity as
\[ SS^\infty(j_! \C_U) = \{ (x, \eta) \in S^*\R^2 \mid x \in \pa U, \eta = d |x| \} = \quad \tikz[baseline=-3pt] {\draw [righthairs] (0,0) circle (.5);}.\]
Here  the Legendrians are represented by co-oriented hypersurfaces in $\R^2$ with hairs indicating the co-orientation. 
\eex

The following Lemma from \cite{KS} is useful in characterising the singular support under $\otimes$ and $\uhom$. 
\bp[Proposition 5.4.14 of \cite{KS}]\label{p:ss-hom-ot}
Let $F$ and $G$ belong to $Sh^b(X)$, then \\
(1) if  $SS^\infty(F) \cap (SS(G)^a)^\infty = \emptyset$, then $SS(F \ot G) \subset SS(F) + SS(G)$. \\
(2) if  $SS^\infty(F) \cap SS^\infty(G) = \emptyset$, then $SS (\uhom(F, G))  \subset SS(G) - SS(F)$. \\
where $(-)^a$ is the fiberwise anti-podal map in $T^*X$ and $\pm$ is the fiberwise sum/substraction in $T^*X$.
\ep
For a version without assuming $SS^\infty(F) \cap SS^\infty(G) = \emptyset$, see Corollary 6.4.5 and 6.2.4 in loc.cit.

\ss{Non-characteristic Deformation Lemma}
Just as in Morse theory, where a level sets $f^{-1}(t)$ of a Morse function $f$ on $M$ has constant diffeomorphism type when $t$ varies in the connected components of the complement of the critical values of $f$, the non-characteristic deformation results for constructible sheaves are about the invariance of the hom-complexes $\hom(F_t, G_t)$ for families of sheaves $\{F_t\}$ and $\{G_t\}$, when $SS^\infty(F)$ and $SS^\infty(G)$ are disjoint.

We first state the version regarding sections of a sheaf over an increasing sequence of open sets. 

\bp[Proposition 2.7.2 in \cite{KS} ] \label{p:ks-nc}
Let $X$ be a real analytic manifold, $F$ a bounded complex of constructible sheaves in $Sh(X)$, and let $\{U_t\}_{t \in \R}$ be a family of open subsets of $X$. We assume the following conditions: \\ 
(1) $U_t = \bigcup_{s < t} U_s$ for all $t \in \R$. \\
(2) For all pairs $(s,t)$ with $s \leq t$, the set $\wb{U_t \RM U_s} \cap \supp(F)$ is compact. \\
(3) Setting $Z_s = \cap_{t>s} \wb{U_t \RM U_s}$ \footnote{See errata at \url{https://webusers.imj-prg.fr/~pierre.schapira/books/Errata.pdf} for the need of closure.}, we have for all pairs $(s,t)$ with $s \leq t$ and all $x \in Z_s \RM U_t$, that 
\[ \uhom(j_{X\RM U_t *} \C_{X\RM U_t},  F)|_x \cong 0. \]
Then we have for all $t \in \R$, the quasi-isomorphism
\[ \Gamma(U, F) \congto \Gamma(U_t; F), \quad \text{ where } U = \bigcup_{s \in \R} U_s. \]
\ep
\brem
The section functor can be viewed as $\Gamma(U_t, F) = \hom(j_{U_t!} \C_{U_t}, F)$.  Hence this is a special case for $\hom(G_t, F_t)$. The advantage for this version is that the results holds for the section over union of the open sets $\{U_s\}$, instead of just between pairs of open sets $U_t, U_s$ for some finite $t,s$. 
\erem

The following two versions of the non-characteristic deformation results are not going to be used in the paper. It is presented here since their conceptual pictures are somewhat clearer. 

\bp[Corollary 2.10, \cite{S}]
\label{p:gamma-inv}
Let $I$ be an open interval of $\R$, let $q: M \times I \to I$ be the projection, and let $\iota_s$ be  the embedding $M \times \{s\} \into M \times I$. Let $F \in Sh(M \times I)$, such that $SS^\infty(F) \cap (T^*_M M \times T^* I)^\infty = \emptyset$ and $q$ is proper on $\supp(F)$. Set $F_s = \iota_s^* F$. Then we have isomorphisms 
\[ \Gamma(M, F_s) \cong \Gamma(M, F_t) \tforall s, t \in I. \]
\ep

Since the hom-complex can be obtained by taking the global section of hom-sheaf, we have a non-characteristic deformation result for $\hom(F_t, G_t)$. First we state a lemma: 
\bl[Petrowsky theorem for Sheaves, Corollary 4.6 \cite{S}] \label{l:pet}
Let $F, G$ be bounded constructibles sheaves in $Sh(X)$. If $SS^\infty(F) \cap SS^\infty(G) = \emptyset$, then the natural morphism 
\[ \uhom(F, \C_X) \ot G \to \uhom(F, G) \]
is an isomorphism. 
\el

\bp\label{p:hom-inv}
Let $I$ be an open interval of $\R$, let $q: M \times I \to I$ be the projection, and let $\iota_s$ be  the embedding $M_s = M \times \{s\} \into M \times I$. Let $F, G \in Sh(M \times I)$, such that \\
(1) $SS^\infty(F)$, $SS^\infty(G)$ and $(T^*_M M \times T^* I)^\infty$ are pairwise disjoint, and\\
(2) $SS^\infty(G_s) \cap SS^\infty(F_s) = \emptyset$, where $F_s = \iota_s^* F, G_s = \iota_s^* G$, and  \\
(3) $q$ is proper on $\supp(F)$ and $\supp(G)$. \\
Then we have isomorphisms 
\[ \hom(F_s, G_s) \cong \hom(F_t, G_t) \tforall s, t \in I. \]
\ep
\bpf
Consider the hom-sheaf $\uhom(F, G)$. We claim that 
\[ SS^\infty( \uhom(F,G) ) \cap (T^*_M M \times T^* I)^\infty = \emptyset. \]
 Indeed, if there exists $((x, t), (p, \tau)) \in T^*(M \times I)$ in $SS( \uhom(F,G)) \cap T^*_M M \times T^* I$ for which $(p, \tau) \neq 0$, then $p=0, \tau \neq 0$ since $T_M^* M = \{(x,0): x \in M\}$.  And there exists $( (x,t), (p_1, \tau_1)) \in SS(F)$ and  $( (x,t), (p_2, \tau_2)) \in SS(G)$, such that $p_1-p_2=0$ and $\tau_1 - \tau_2=\tau$ (Proposition \ref{p:ss-hom-ot}). By condition (1), $p_1 \neq 0$ and $p_2 \neq 0$. By condition (2), if $p_1,p_2$ are non-zero then $p_1 \neq p_2$. Hence it is impossible to have $p_1 - p_2=0$, and the claim is proven. 

Hence by Proposition \ref{p:gamma-inv}, we have
\[ \Gamma(M, \uhom(F,G)|_s) \cong \Gamma(M, \uhom(F,G)|_t) \tforall t, s\in I. \]

Finally, we claim that
\[ \hom_M(\iota_s^* F, \iota_s^* G) \cong  \Gamma(M, \uhom(F,G)|_s) \tforall  s\in I. \]
Let $\C_t: = \iota_{t*} \iota_t^* \C_{M \times I}$ be the constant sheaf support on the slice $M \times \{t\}$, for any $t \in I$. Then $SS^\infty(\C_t)$ is disjoint from $SS^\infty(F), SS^\infty(G)$ and $SS^\infty(\uhom(F,G))$  by assumption and the first claim. Hence we have 
\[ \C_s \ot F \cong \uhom( \uhom(\C_s, \C_{M \times I}), \C_{M \times I}) \ot F \cong  \uhom( \C_s[-1], \C_{M \times I}) \ot F \congto \uhom(\C_s[-1], F) \]
where Petrowsky Theorem is applied in the last isomorphism. The same is true by replacing $F$ with $G$ and $\uhom(F,G)$. Then
\bea 
 && \hom_M(\iota_s^* F, \iota_s^* G) \cong \hom_{M \times I} (F, \iota_{s *} \iota_s^* G) \cong \hom_{M \times I} (F, \C_s \ot G) \\
 &\cong & \hom_{M \times I} (F, \uhom(\C_s[-1], G)) \cong  \hom_{M \times I} (\C_s[-1], \uhom(F, G)) \\
 &\cong& \Gamma(M\times I, \uhom(\C_s[-1], \C_{M\times I}) \ot \uhom(F, G)) =  \Gamma(M\times I, \C_s \ot \uhom(F, G)) \\
 & \cong &\Gamma(M, \uhom(F, G)|_s)
\eea
This finishes the proof of the proposition. 
\epf

\section{Toric Geometry}
An $n$-dimensional  smooth projective complex manifold $X$ is toric  if there is a holomorphic $(\C^*)^n$-action with an open dense orbit $X^o$ on which $(\C^*)^n$ acts freely. The complement of the open orbit $D=X \RM X^o$ is a simple normal crossing divisor with irreducible torus-invariant components. 

We first review the standard setup and notation for the combinatorial data used for defining a toric variety. Then we explain the relationship between equivariant line bundles, toric divisors, and twisted polytopes (as a collection of labeled vertices). Finally, we review `twisted polytope sheaves', the corresponding constructible sheaves for equivariant line bundles under the equivariant CCC. 

\ss{Fan Data} The data of a toric manifold can be expressed combinatorically using a fan. Let $N \cong \Z^n$ be a rank $n$ lattice, with $N_\R = N \ot_\Z \R$.  Let $M = \Hom(N, \Z)$ be the dual lattice and $M_\R = M \ot_\Z \R$ be the dual vector space. Let $\la -, - \ra: M_\R \times N_\R \to \R$ be the dual pairing.  Let $T_M = M_\R / M$ be a real $n$-dimensional torus, and $\pi: M_\R \to T_M$ be the quotient map.  We recall the following definitions. 
\begin{enumerate}
\item  A convex polyhedral cone $\sigma \subset N_\R$ is a set of the form $\sigma = \cone(S) = \{\sum_{u \in S} \lambda_u u \mid \lambda_u \geq 0 \}$, where the cone generator $S \subset N_\R$ is a finite subset. A cone $\sigma$ is {\em rational} if there is a generator $S$ for $\sigma$ such that $S \subset N$. A cone is {\em strongly convex} if it does not contain any non-trivial linear subspace of $N_\R$. 
\item Let $\sigma \in \Sigma$ be a cone, we define the dual (closed) cone $\sigma^\vee$ as
\[ \sigma^\vee: = \{ x \in M_\R \mid \la x, y \ra \geq 0, \forall y \in \sigma \}. \]
We also define $\sigma^\perp = \{ x \in M \mid \la x, y \ra = 0, \forall y \in \sigma \} \subset M_\R$, and $\sigma^o$ (resp. $(\sigma^\vee)^o$) as the relative interior of $\sigma$ (resp. $\sigma^\vee$).
\item A face of a cone $\sigma$ is the subset $H_m \cap \sigma$ for some $m \in \sigma^\vee$ and $H_m = m^\perp$. We use $\sigma(r)$ to denote the collection of $r$-dimensional faces of $\sigma$. 
\item  A fan $\Sigma$ in $N_\R$ is a finite collection of strongly convex rational polyhedral cones $\sigma \subset N_\R$, such that (a) if $\sigma \in \Sigma$ then any face of $\sigma$ is in $\Sigma$, and (b) if $\sigma_1,\sigma_2$ are cones in $\Sigma$ then $\sigma_1 \cap \sigma_2$ is a face in both $\sigma_1$ and $\sigma_2$. We use $\Sigma(r)$ to denote the collection of $r$-dimensional cones in $\Sigma$. 
\item A fan $\Sigma$ in $N_\R$ is {\em complete}, if its support $|\Sigma| := \cup_{\sigma \in \Sigma} \sigma$ is the entire $N_\R$.  A complete fan $\Sigma$ is {\em smooth}, if each maximal cone $\sigma \in \Sigma(n)$ is generated by a lattice basis of $N$. 
\item A smooth complete fan $\Sigma$ is {\em projective}, if there exists a convex piecewise linear function $\varphi: N_\R \to \R$, such that the maximal linearity domains of $\varphi$ are the maximal cones of $\Sigma$. (cf. Proposition \ref{p:divisor})
\end{enumerate}
See Example \ref{ex:p2} for a fan of $\P^2$. 

\noindent {\bf Assumption:} We will always assume $\Sigma$ to be a smooth projective fan.

The affine toric variety $X_\sigma$ is then defined by
\[ X_\sigma = \Spec (\C[\sigma^\vee \cap M ]) \]
where $\C[\sigma^\vee \cap M ])$ is the group ring of the abelian semi-group $\sigma^\vee \cap M$. 
If $\tau \subset \sigma$ is a face of $\sigma$, then $\sigma^\vee \subset \tau^\vee$, hence $\C[\sigma^\vee \cap M ]) \into \C[\tau^\vee \cap M ])$, and $X_\tau \into X_\sigma$ is an open inclusion. We may equip $\Sigma$ with a partial ordering, for $\tau, \sigma \in \Sigma$,  $\tau \leq \sigma \iff \tau \subset \sigma$. Then $X_\Sigma$ can be glued together from affine open pieces $X_\sigma$, as a colimit of schemes
\[ X_\Sigma = \colim_{\sigma \in \Sigma} X_\sigma. \]



\ss{Toric Divisors, Support Functions and Twisted Polytopes}  
For each ray $\rho \in \Sigma(1)$,  let $v_\rho \in \rho \cap N$ be a minimal ray generator, $\lambda^{v_\rho}: \C^* \to N \ot_\Z \C^*$ as the one-parameter subgroup, and $D_\rho = \overline{ \{\lim_{t\to 0} \lambda^{v_\rho}(t) \cdot x \mid x \in X^o \} }$ the torus-invariant divisor, or, a {\em toric divisor}. We write $\sum_\rho$ for a summation over the rays $\rho \in \Sigma(1)$ when there is no danger of confusion. 

Let 
$ D = \sum_\rho a_\rho D_\rho$
be a toric $\R$-divisor on $X_\Sigma$, $a_\rho \in \R$. If $a_\rho \in \Z$ for all $\rho$, then $D$ is an integral toric divisor, or toric $\Z$-divisor.  There are two equivalent ways to describe a toric divisor, either using a support function $\varphi_D$ on $N_\R$, or a twisted polytope $\chi_D$ on $M_\R$.

\bd[Support function]
A {\em support function} for $\Sigma$ is a continuous piecewise linear function $\varphi: N_\R \to \R$, such that for each maximal cone $\sigma \in \Sigma(n)$, the restriction $\varphi|_\sigma$ is linear.  
\begin{itemize}
\item A support function $\varphi$ is integral if it sends $N$ to $\Z$.  
\item A support function $\varphi$ is  {\em convex}, if for any $x,y \in N_\R$
\[\varphi(tx + (1-t)y) \leq t \varphi(x) + (1-t) \varphi(y). \]
\item Furthermore, we say $\varphi$ is {\em strictly convex}, if the strict inequality holds whenever $x,y$ is not contained in the same cone. 
\end{itemize}
\ed

\bd[Twisted polytope] \label{d:tp}
A {\em twisted polytope} for $\Sigma$ is an assignment of element in $M_\R$ to  top-dimensional cones in $\Sigma$, 
\[ \chi: \Sigma(n) \to M_\R, \quad \sigma \mapsto \chi_\sigma, \] 
such that if $\sigma, \tau \in \Sigma(n)$
then $\la \chi_\sigma,  \cdot \ra = \la   \chi_\tau,  \cdot \ra $ on  $\sigma \cap \tau$. 
\begin{itemize}
\item A twisted polytope $\chi$ is  integral   if the function $\chi_\sigma  \in M$ for all $\sigma \in \Sigma(n)$.  
\item If $\chi$ is a twisted polytope, then for any cone $\sigma \in \Sigma$,  we define $\chi_\sigma \in M_\R/\sigma^\perp$ by
\[ \chi_\sigma = \text{ Affine Hull} (\{ \chi_\tau \mid \text{ $\tau$ is a maximal cone containing $\sigma$} \}) \subset M_\R. \]
\item For any $x \in M_\R$, let $\chi+x$ denote the translated twisted polytope that sends $\sigma \mapsto \chi_\sigma + x$ for any $\sigma \in \Sigma(n)$. 
\end{itemize}
\ed
\brem
The data for a twisted polytope is a collection of the vertices, labelled by $\Sigma(n)$.  
\erem

\bp
Let $\Sigma$ be an $n$-dimensional smooth projective fan. Then we have a canonical equivalences among the following three types of objects. 
(1) A toric $\R$-divisor $D= \sum_\rho a_\rho D_\rho$, $a_\rho \in \R$. \\
(2) A twisted polytope, $\chi: \Sigma(n) \to M_\R$. \\
(3) A support function, $\varphi: N_\R \to \R$. \\
In particular, integral toric divisors corresponds to integral twisted polytopes and integral support functions. 
\ep
\bpf
$(2) \LRA (3)$. Given $\chi$, we may define $\varphi$ by $\varphi(x) = \la \chi_\sigma, x\ra$ if $x \in \sigma$ for some maximal cone $\sigma \in \Sigma(n)$. This is well-defined since if $x \in \sigma \cap \tau$, then $\la \chi_\sigma,x \ra = \la  \chi_\tau,x \ra$. Conversely, given $\varphi$, then for each maximal cone $\sigma$, the linear function $\varphi|_{\sigma}$ determines an element in $M_\R$, denoted by $\chi_\sigma$. The continuity of $\varphi$ ensures $\la \chi_\sigma,  \cdot \ra = \la   \chi_\tau,  \cdot \ra $ on  $\sigma \cap \tau$.

$(1) \LRA (3)$. Given a toric $\R$-divisor $D = \sum_{\rho} a_\rho D_\rho$, for each $\rho \in \Sigma(1)$, we define 
$\varphi|_\rho: \rho \to \R$ by $v_\rho \mapsto a_\rho$. Since the cones of $\Sigma$ are simplicial, there is a unique piecewise linear extension of $\varphi$ to $N_\R$ that is linear in each cone of $\Sigma$. Conversely, given $\varphi$, let $a_\rho = \varphi(v_\rho)$ for each $\rho \in \Sigma(1)$. 

The claim on integrality is straightforward to verify. This finishes the proof of the Proposition. 
\epf

\brem
(1) If $D$ is a toric $\R$-divisor, we let $\chi_D$ and $\varphi_D$ be the corresponding twisted polytope and support function. 
(2) If $\chi$ is a twisted polytope for $\Sigma$, and $\varphi$ is the corresponding support function, then $\chi_\sigma  \in M_\R / \sigma^\perp$ corresponds to the linear function $\varphi|_\sigma: \sigma \to \R$. 
\erem

\bp \label{p:divisor}
Let $X_\Sigma$ be a smooth complete toric variety. Let $D = \sum_{\rho} a_\rho D_\rho$ be an integral toric divisor. Then \\
(1) $D$ is base-point free if and only if $\varphi_D$ is convex.  \\
(2) $D$ is ample if and only if $\varphi_D$ is strictly convex. 
\ep
\bpf
\cite{CLS}, Chapter 4 and 6. 
\epf

\bd
If $D=\sum_{\rho} a_\rho D_\rho$ is ample, we define the open  convex polytope $\Delta_D$ as the interior of the convex hull of $\{\chi_\sigma \mid \sigma \in \Sigma(n)\}$. Equivalently, we have
\[ \Delta_D = \{ x \in M_\R \mid \la x, v_\rho \ra < a_\rho, \text{ for all } \rho \in \Sigma(1) \}. \] 
\ed

\ss{Constructible Sheaves and Twisted Polytope Sheaves}
Let $M, N$ be dual rank-$n$ lattices, and $\Sigma$ a smooth complete fan in $N_\R$. We define the conical Lagrangians $\Lambda_\Sigma$ in $T^* M_\R$ as \footnote{Our definition differs in sign convension from that in \cite{FLTZ}. If we change $\Sigma$ to $-\Sigma$ in this paper, then the definition agrees.}
\be
\label{d:lam-r} \Lambda_{\Sigma} = \bigcup_{\sigma \in \Sigma} (\sigma^\perp + M)  \times \sigma \subset M_\R \times N_\R = T^* M_\R.
\ee
We denote the push-forward of $\Lambda_\Sigma$ to $T^* T_M$ by $\wb \Lambda_\Sigma$, or directly we have
\be
\label{d:lam-t} \wb \Lambda_{\Sigma} = \bigcup_{\sigma \in \Sigma} (\sigma^\perp / \sigma^\perp \cap M)  \times \sigma \subset T_M \times N_\R = T^* T_M.
\ee

\bd[Standard Shard Sheaves]
For any cone $\sigma \in \Sigma$, $c \in M_\R / \sigma^\perp$,  we define the closed subset $Q({\sigma, c}) \subset M_\R$ and the standard sheaf $P(\sigma,c)$ as
\[  Q(\sigma, c):=c + \sigma^\vee \subset M_\R,  \quad P(\sigma, c) := j_{ Q(\sigma, c) *} \C_{ Q(\sigma, c)}.  \]
\ed

\bd[Twisted Polytope Sheaves on $M_\R$]
\label{d:pshv}
Let $\chi$ be a twisted polytope for $\Sigma$, let $D$ be the corresponding toric $\R$-divisor. The {\em twisted polytope sheaf $ P(\chi)$ on $M_\R$ } is defined by the following chain complex of sheaves, with $\C_{M_\R}$ at degree $-n$, 
\[ P(\chi) := \left( \C_{M_\R} \xto{d_1} \bigoplus_{\sigma_1 \in \Sigma(1)} P(\sigma_1, \chi_{\sigma_1})  \xto{d_2}  \bigoplus_{\sigma_2 \in \Sigma(2)} P(\sigma_2, \chi_{\sigma_2})  \xto{d_3} \cdots \xto{d_n} \bigoplus_{\sigma_n \in \Sigma(n)} P(\sigma_n, \chi_{\sigma_n}) \right) \] 
where   $d_k$ for $k=1,\cdots, n$ is given in the following way:
\[ d_k = \sum_{\sigma_{k-1} \subset \sigma_{k}} sgn(\sigma_{k-1}, \sigma_{k}) \rho_{\sigma_{k} \to \sigma_{k+1}} \]
where the sum is over  $\sigma_{k-1} \in \Sigma(k-1)$,  $\sigma_{k} \in \Sigma(k)$, and
\[ \rho_{\sigma_{k-1} \to \sigma_{k}}: P(\sigma_{k-1}, \chi_{\sigma_{k-1}}) \to  P(\sigma_{k}, \chi_{\sigma_{k}}) \]
is the canonical restriction, and the sign $sgn(\sigma_{k-1}, \sigma_k) = \pm 1$ is chosen such that $d^2=0$ (see the following remark). 

If $D$ is any toric $\R$-divisor, $\chi=\chi_D$ the corresponding twisted polytope,  we sometimes write $P(D)$ for $P(\chi_D)$. 
\ed
\brem
We can be more concrete about the sign choices $sgn(\sigma_{k-1}, \sigma_k)$. One way is to fix a linear ordering of the rays $\Sigma(1)$, then a $k$-dimensional simplicial cone $\sigma_k$ can be identified with the ordered set $\sigma_k(1) = \{\rho_1<\rho_2<\cdots<\rho_k\}$. If $\sigma_{k-1} = \sigma_k - \{\rho_j\}$, then we set $sgn(\sigma_{k-1}, \sigma_k)=(-1)^{j-1}$. Another way is to fix the orientations of all cones in $\Sigma$ once and for all, and $sgn(\sigma_{k-1}, \sigma_k)=\pm 1$ depending on if $\sigma_{k-1}$ agrees with the induced boundary orientation of $\sigma_k$. 
\erem

\bex
Consider the example of $\P^1$, where $\Sigma(1)=\{\R v_1, \R v_2\}$, where $v_1 = 1$ and $v_2 = -1$. We still need to fix the `offset parameters' $\chi_i$ for each $v_i$. We consider the following three cases
\begin{enumerate}
\item $\chi_1=-1, \chi_2 = 1$, then 
\[ P(\chi) \cong ( \C_{\R} \to \C_{[-1, \infty)} \oplus \C_{(-\infty, 1]} ) \cong \C_{[-1,1]} \]
\item $\chi_1=0, \chi_2 = 0$, then 
\[ P(\chi) \cong ( \C_{\R} \to \C_{[0, \infty)} \oplus \C_{(-\infty, 0]} ) \cong \C_{\{0\}} \]
\item $\chi_1=1, \chi_2 = -1$, then 
\[ P(\chi) \cong ( \C_{\R} \to \C_{[1, \infty)} \oplus \C_{(-\infty, -1]} ) \cong \C[1]_{(-1,1)} \]
\end{enumerate}
\begin{figure}
%
%
%
\includegraphics[width=\textwidth]{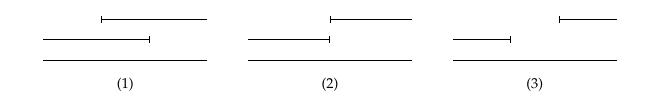};
\caption{Twisted Polytope Sheaves for $\P^1$ \label{f:p1}}. 
\end{figure}
where we briefly abuse notation and denote by $\C_A$ the constant sheaf supported on the subset $A$. The supports of the standard sheaves in the chain complexes also shown in Figure \ref{f:p1}. 
\eex

Since $M_\R$ is a vector space, we have the addition operation $v: M_\R \times M_\R \to M_\R$. 
The addition operation induces the convolution product $\star$ for sheaves $Sh(M_\R)$
\[ F_1 \star F_2 : = v_! ( F_1 \boxtimes F_2). \]

We have the following properties of twisted polytope sheaves. 
\bp\label{p:tw-prop}
Let $\Sigma$ be a smooth projective fan, $D = \sum_\rho a_\rho D_\rho$ a toric $\R$-divisor, and $P(D)$ the twisted polytope sheaves on $M_\R$. Then
\begin{enumerate}
\item If $D$ is  integral, then there is a unique up to isomorphism equivariant line bundle $\ocal_X(D)$ on $X_\Sigma$, and
\[ \kappa(\ocal_X(D)) = P(D).\] 
In particular $P(0) =j_{\{0\} *} \C_{\{0\}}$ is the skyscraper sheaf at point $0$. 
\item If $D$ is an ample divisor, then $P(D)$ is a costandard sheaf supported on a simplicial convex polytope, with each facet corresponding to a ray $\rho \in \Sigma(1)$, and each vertex corresponding to a maximal cone $\sigma \in \Sigma(n)$. 
\item If $-D$ is an ample divisor, then $P(D)$ is a standard sheaf supported on a simplicial convex polytope. $P(D)  = a_* \D (P(-D))$, where $a: M_\R \to M_\R$ sends $x \mapsto -x$. 
\item $P(D)$ has compact support in $M_\R$. For any $x \in M_\R$, the stalk $P(D)_x$ has cohomology in degrees between $-n$ and $0$. 
\item If $D = D_1 + D_2$, then $P(D) = P(D_1) \star P(D_2)$, where $\star$ is the convolution product on $M_\R$.
\item Let $\chi$ be any twisted polytope, then $(-) \star P(\chi): Sh(M_\R, \Lambda_\Sigma) \to Sh(M_\R, \Lambda_\Sigma)$ is an equivalence of cateogry. The functor $(-) \star P(\chi)$  has an inverse  $(-) \star P(-\chi)$ 
\end{enumerate}
\ep
\bpf
The results are given in \cite{FLTZ, Tr}, with straightforward adaptations from integer to real coefficients. 
\epf

\bl
Let $D=\sum_\rho a_\rho D_\rho$.  If $a_\rho$ is not an integer for any $\rho \in \Sigma(1)$, then $SS^\infty(P(D)) \cap \Lambda_\Sigma^\infty = \emptyset$. 
\el
\bpf
From the chain complex definition for $P(\chi)$, we have
\[ SS( P(\chi)) \subset \bigcup_{\sigma \in \Sigma} SS(P(\sigma, \chi_\sigma)) =  \bigcup_{\sigma \in \Sigma} (\chi_\sigma + \sigma^\perp) \times \sigma. \]
If $(x, p) \in SS( P(\chi)) \cap \Lambda_\Sigma$ and $p \neq 0$, then there are non-zero cones $\sigma, \tau \in \Sigma$, such that
\[ (x,p) \in \left( (\chi_\sigma + \sigma^\perp) \times \sigma \right)\bigcap \left(  (M + \tau^\perp) \times \tau\right). \]
Hence $p \in \sigma \cap \tau$.Thus $\sigma \cap \tau$ contains at least a ray $\rho \in \Sigma(1)$, otherwise $p=0$. Consider $\la x, v_\rho \ra$.  Since $x \in \chi_\sigma + \sigma^\perp$, we have
$ \la x, v_\rho \ra = a_\rho$. On the other hand, $x \in M + \tau^\perp$,  hence $\la x, v_\rho \ra \in \Z$. This contradicts with $a_\rho \notin \Z$ for any $\rho \in \Sigma(1)$. Thus the lemma is proven. 
\epf

\bd[Twisted Polytope Sheaves on $T_M$]
Let $\chi$ be a twisted polytope for $\Sigma$, $P(\chi)$ the twisted polytope sheaf for $\chi$ on $M_\R$, then the twisted polytope sheaf for $\chi$ on $T_M$ is
\[ \wb P(\chi) := \pi_*  P(\chi) = \pi_! P(\chi). \]
where $\pi_* = \pi_!$ since $\pi$ is proper on $\supp P(\chi)$.  
\ed
\brem \label{rem:lift}
For any lattice point $x \in M$, the shifted polytope $\chi + x$ defines the same twisted polytope sheaf,  $\wb P(\chi) = \wb P(\chi+x)$, since $\pi = \pi \circ ( \cdot + x): M_\R \to T_M$. 
\erem
The following is stated in \cite{Tr}.
\bp
If $D$ is an integral twisted polytope, then the non-equivariant CCC functor $\wb \kappa$ sends $\ocal_X(D)$ to $\wb P(D)$. 
\ep

\bex \label{ex:p2}
Consider the following two dimensional fan $\Sigma$, with ray generators $v_1 = (1,0), v_2 = (0,1), v_3 = (-1,-1)$. Let $D = D_1 + D_2 + D_3$ where $D_i$ is the toric divisor for the ray $v_i$, then $\varphi_D$ is a strictly positive function on $N_\R$, such that $\varphi_D(v_i) =1$. The vertices for the twisted polytope $\chi_D$ are $(1,1), (-2,1), (1,-2)$ in $M_\R$. The twisted polytope sheaf $ P(\chi)$ is the costandard sheaf supposed on the interior of the shaded region. 
The blue hairs indicate the singular support $SS^\infty(P(\chi))$ at infinity. 
%
%

\begin{figure}[h]
\includegraphics[width=\textwidth]{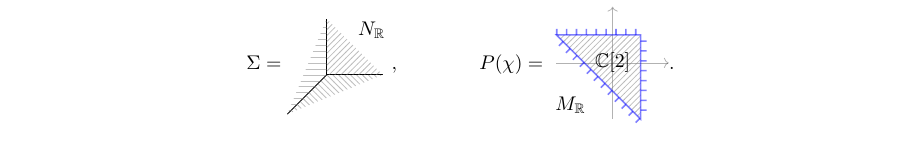};
\end{figure}


\eex

\section{Proof of Main Theorem}
Let $\pcal$ be the full dg subcategory of $Sh(M_\R, \Lambda_\Sigma)$ spanned by the {\em integral} twisted polytope sheaves on $M_\R$, and let $\wb \pcal$ be the full dg subcategory of $Sh(T_M, \wb \Lambda_\Sigma)$ spanned by the {\em integral} twisted polytope sheaves on $T_M$.

The follow proposition is the heart of this paper. We first define the probe sheaves for the stalk over $x \in M_\R$ and $\theta \in T_M$. 
\bd \label{d:tw-div}
For any $x  \in M_\R$, let the integral toric divisor $D_{[x]}$ be defined by
\[  D_{[x]} : = \sum_{\rho} (\lfloor \la x, v_\rho \ra \rfloor + 1) D_\rho, \quad P_{[x]} := P(D_{[x]}).\]
For any $\theta \in T_M$, we may fix any lift  $x$ of $\theta$ in $M_\R$, then define
\[  \wb P_{[\theta]} := \pi_* P_{[x]}. \]
Since different lifts of $x$ differ by an element in $M$, hence the push-forward is independent of the choice of the lift. (cf. Remark \ref{rem:lift}. )
\ed

\bp\label{p:stalk}
For any point $\theta \in T_M$, there is a unique twisted polytope sheaf $\wb P_{[\theta]}$ on $T_M$, such that for any sheaf  $\wb F \in Sh(T_M, \wb \Lambda_\Sigma)$, the stalk at $\theta$ can be computed by
\[ \wb F_\theta \cong \hom(\wb P_{[\theta]}[-n], \wb F). \]
\ep

\bpf
Fix any $x \in \pi^{-1}(\theta)$. 
For any sheaf  $\wb F \in Sh(T_M, \wb \Lambda_\Sigma)$, let $F = \pi^{-1} \wb F= \pi^! \wb F$. Then we have canonical isomorphisms 
\[ F_x \cong \wb F_\theta, \]
and
\[   \hom(\wb P_{[\theta]}[-n], \wb F) \cong  \hom(\pi_! P_{[x]}[-n], \wb F) \cong \hom( P_{[x]}[-n], \pi^! \wb F) =  \hom( P_{[x]}[-n], F). \]
Hence it suffices to prove that for any fixed $x \in M_\R$, we have
\be \tag{*} F_{x} \cong \hom( P_{[x]}[-n], F). \ee

Since $\Sigma$ is smooth projective, there exists an integral ample toric divisor
\[ A = \sum_{\rho} a_\rho D_\rho, \quad a_\rho \in \Z_{>0}. \]
Then the twisted polytope sheaf $P(A)$ is supported on $\Delta_A$, with stalk $\C[n]$. Since $a_\rho >0$, we have $0 \in \Delta_A$. 

Fix $\epsilon_0>0$ small enough, depending only on $x$ and $\Lambda_\Sigma$,  such that for any $0<\epsilon \leq \epsilon_0$, 
\[ F_x \cong \Gamma( \Delta_{D_{(x)} + \epsilon A}, F) \]
where $ \Delta_{D_{(x)} + \epsilon A}, = x + \epsilon \Delta_A$ is a shifted open convex polytope around $x$. This is possible since $F$ is a polyhedral constructible sheaf, and $\Delta_A$ is a convex set. In particular, we may shrink $\epsilon_0$ and further assume that 
\[ \epsilon_0 a_\rho + \la x, v_\rho \ra  < \lfloor \la x, v_\rho \ra \rfloor + 1, \tforall \rho \in \Sigma(1). \]

Fix $R > 0$ a large enough integer, such that $D_{[x]} + R A$ is an ample integral toric divisor. Let $\Delta_{D_{[x]} + R A}$ be the corresponding open convex polytope. 
 
For any $s \in [0,1]$, we define a 1-parameter family of ample divisors $D_s$, interpolating between $x +(R+\epsilon_0) A$ and $D_{[x]} + RA$. 
\[ D_s = \sum_{\rho} a_{\rho,s} D_\rho, \quad a_{\rho,s} = (1-s) (\la x, v_\rho \ra + (R+\epsilon_0) a_\rho) + s ( \lfloor \la x, v_\rho \ra \rfloor + 1 + R a_\rho),\]
and let 
\[ \Delta_s := \Delta_{D_s},\quad \text{ and } \quad P_s := P(D_s). \] 
Since $\lfloor \la x, v_\rho \ra \rfloor + 1 > \la x, v_\rho \ra + \epsilon_0 a_\rho > \la x, v_\rho \ra$, and there is no integer in the open interval $(\la x, v_\rho \ra, \lfloor \la x, v_\rho \ra \rfloor + 1)$, hence for any $s \in (0,1)$, 
\[ \la x, v_\rho \ra + \epsilon_0 a_\rho< a_{\rho,s}  - R a_{\rho} < \lfloor \la x, v_\rho \ra \rfloor + 1, \text{ and } a_{\rho, s} \notin \Z. \]
Thus from Lemma \ref{p:tw-prop} , 
\[ SS^\infty( P_s) \cap \Lambda_\Sigma^\infty = \emptyset \tforall s \in (0,1). \]

Apply the non-characteristic deformation result in Proposition \ref{p:ks-nc}, let 
\[ U_t = \bcs \Delta_0, & t \leq 0 \\
\Delta_{t/(1+t)}, & t > 0
\ecs
\]
we have for any sheaf $G \in Sh(M_\R, \Lambda_\Sigma)$, 
\[ \Gamma( \Delta_1, G) = \Gamma(\cup_{s \in (0,1)} \Delta_s, G) \cong \Gamma(\Delta_t, G) \tforall t \in [0,1). \]
Since $D_s$ are ample divisor for all $s \in [0,1]$ (since ample cone is convex),  and $P(D_s) = j_{\Delta_s !} \omega_{\Delta_s} \cong j_{\Delta_s !} \C_{\Delta_s}[n]$, we have
\[ \Gamma(\Delta_s, G) = \hom(j_{\Delta_s !} \C_{\Delta_s}, G) \cong \hom(P(D_s)[-n], G),  \tforall s \in [0,1]. \]

Finally, we use convolution $\star$ is an equivalence of category on $Sh(M_\R)$ to get
\bea 
F_x \cong \Gamma( x + \epsilon \Delta_A, F) &\cong& \hom( j_{x + \epsilon \Delta_A !} \C_{x + \epsilon \Delta_A}, F)  \\
&\cong & \hom( j_{x + \epsilon \Delta_A !} \C_{x + \epsilon \Delta_A} \star P(RA), F \star P(RA)) \\
& \cong & \hom(P(D_0)[-n], F \star P(RA)) \\
& \cong & \hom(P(D_1)[-n], F \star P(RA)) \\
& \cong & \hom(P(D_1-RA)[-n], F) \\
& \cong & \hom(P(D_{[x]})[-n], F).\\
\eea
This finishes the proof of the Proposition. 
\epf

Now we prove the main theorem stated in the introduction section. 
\bpf[Proof of the Main Theorem]
First we claim that there exists a semi-orthogonal expansion  
\[ Sh(T_M, \wb \Lambda_\Sigma) \cong \la  \la \wb \pcal \ra^\perp  ,  \la \wb \pcal \ra   \ra. \]
Since $\wb \kappa$ is an quasi-embedding, hence 
\[ Coh(X_\Sigma) \congto \wb \kappa( Coh(X_\Sigma)). \]
Since line bundles generates $Coh(X_\Sigma)$, hence
\[ \wb \kappa( Coh(X_\Sigma)) \cong \wb \kappa( \la \{ \lcal : \text{ line bundles} \} \ra ) \cong \la \wb \pcal \ra. \]
Kawamata proved that $Coh(X_\Sigma)$  admits an exceptional collection \cite{Ka}, hence $\la \wb \pcal \ra$ also admits an exceptional collection. By Proposition 2.6 and Corollary 2.10 in \cite{BK}, $\la \wb \pcal \ra$ is saturated and is left and right admissible. In other words, the semi-orthogonal decomposition in the claim exists. 

From Proposition \ref{p:stalk}, we have $\la \wb \pcal \ra^\perp=0$. Hence $Sh(T_M, \wb \Lambda_\Sigma) \cong \la \wb \pcal \ra \cong Coh(X_\Sigma)$.
\epf

\ifpaper

\bex
We consider two toric surfaces, with some twisted polytopes $P_{(x)}$ shown in Figure \ref{fig:tw}. 
\begin{enumerate}
\item Let $X_\Sigma = \P^2$.  The red, blue and yellow twisted polytopes  are  $\ocal(3), \ocal(2), \ocal(1)$ respectively. 
\item Let $X_\Sigma = \F_3$, with ray generators $(1,0), (0,1), (-1,-3), (-1,0)$. This is a smooth non-Fano projective toric surface. The red polytope corresponds to the anti-canonical bundle, with all $a_\rho =1$. Indeed, it is non-Fano since the anti-canonical bundle is twisted. The yellow polytope is ample. 
\end{enumerate}

\begin{figure}[h]
\includegraphics[width=\textwidth]{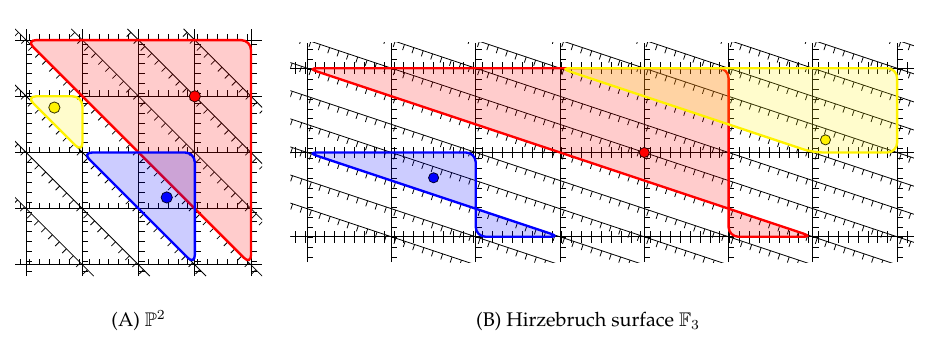}
\caption{\label{fig:tw} Various probe sheaves $P_{[x]}$ (shown as colored twisted polytopes) on $M_\R$ for different $x$ (shown as solid dots).}
\end{figure}
\eex

\fi


\begin{thebibliography}{HHHH}

\bibitem[BK]{BK} A. I. Bondal, M. M. Kapranov. Representable functors, Serre functors, and mutations. MATH USSR IZV, 1990, 35 (3), 519–541

\bibitem[CLS]{CLS} David A. Cox, John B. Little, and Henry K. Schenck. Toric Varieties. American Mathematical Soc., Jan 1, 2011 

\bibitem[CS]{CS} Alastair Craw, Gregory G. Smith. Projective toric varieties as fine moduli spaces of quiver representations. American Journal of Mathematics 130 (2008) 1509-1534



\bibitem[Dr]{Dr} V. Drinfeld. DG Quotients of DG categories.  Journal of Algebra Volume 272, Issue 2, 15 February 2004, Pages 643-691

\bibitem[E]{E} Alexander Efimov. Maximal lengths of exceptional collections of line bundles. Journal of the London Mathematical Society 90(2) · October 2010

\bibitem[FLTZ]{FLTZ} Bohan Fang, Chiu-Chu Melissa Liu, David Treumann, Eric Zaslow. A categorification of Morelli's theorem. 	Invent. Math. 186 (2011), no.1, 79-114





\bibitem[HP]{HP} Lutz Hille and Markus Perling, A counterexample to King’s conjecture, Compos. Math. 142 (2006), no. 6, 1507–1521

\bibitem[Ka]{Ka} Yujiro Kawamata. Derived categories of toric varieties. Michigan Math. J.
Volume 54, Issue 3 (2006), 517-536.

\bibitem[Ke]{Ke} Bernhard Keller. On differential graded categories. arXiv:math/0601185

\bibitem[KS]{KS} Masaki Kashiwara, Pierre Schapira, Sheaves on Manifolds. 

\bibitem[KT]{KT} Yael Karshon and Susan Tolman, The moment map and line bundles over presymplectic toric manifolds. J. Differential Geom. Volume 38, Number 3 (1993), 465-484.

\bibitem[Ku1]{Ku1} Tatsuki Kuwagaki, The nonequivariant coherent-constructible correspondence for toric surfaces, arXiv preprint arXiv:1507.05393, to appear in Journal of Differential Geometry (2015).

\bibitem[Ku2]{Ku2} ---------, The nonequivariant coherent-constructible correspondence for toric stacks, arXiv 1610.03214

\bibitem[N]{N} David Nadler. Microlocal branes are constructible sheaves, Selecta Math. 15 (2009), no. 4, 563--619.



\bibitem[NZ]{NZ} David Nadler, Eric Zaslow. Constructible Sheaves and the Fukaya Category. J. Amer. Math. Soc. 22 (2009), 233-286 

\bibitem[SS]{SS} Sarah Scherotzke and Nicol`o Sibilla, The non-equivariant coherent-constructible correspondence and a conjecture of King, Selecta Math. (N.S.) 22 (2016), no. 1, 389–416. MR 3437841

\bibitem[S]{S} Pierre Schapira. A short review on microlocal sheaf theory. \href{https://webusers.imj-prg.fr/~pierre.schapira/lectnotes/MuShv.pdf}{Listed on personal website}. 

\bibitem[STW]{STW} Vivek Shende, David Treumann, Harold Williams. On the combinatorics of exact Lagrangian surfaces. arXiv:1603.07449

\bibitem[Tr]{Tr} David Treumann, Remarks on the nonequivariant coherent-constructible correspondence for toric varieties, arXiv preprint arXiv:1006.5756 (2010).


\end{thebibliography}
\end{document}